\documentclass{article}
\usepackage{graphicx} 
\usepackage{amsfonts}
\usepackage{amsmath}
\usepackage{amssymb}
\usepackage{amsthm}
\usepackage{mathtools}
\usepackage{hyperref}
\usepackage[frozencache=true,cachedir=minted-cache]{minted}
\usepackage{cleveref}
\usepackage[utf8]{inputenc}
\usepackage[mathletters]{ucs}
\usepackage[normalem]{ulem}

\usepackage{newunicodechar}
\newunicodechar{μ}{\textmu}
\newunicodechar{ε}{\textepsilon}
\newunicodechar{ω}{\textomega}
\newunicodechar{⁴}{\ensuremath{^4}}
\newunicodechar{β}{\textbeta}
\usepackage{textgreek}

\usepackage[usenames,dvipsnames,x11names]{xcolor}
\definecolor{burntumber}{rgb}{0.54, 0.2, 0.14}

\usepackage[format=hang,font=small,labelfont={bf,sl},textfont=sl,width=0.95\textwidth]{caption}
\usepackage{subcaption}
\usepackage[body={155mm,226mm},headsep=20mm,marginparsep=4mm,hmarginratio=1:1,bottom=1.5cm]{geometry}

\newtheorem{theorem}{Theorem}[section]

\theoremstyle{definition}

\newtheorem{definition}[theorem]{Definition}

\newcommand{\bol}{\boldsymbol}

\newcommand{\nef}{\boldsymbol{f}}
\newcommand{\neu}{\boldsymbol{u}}
\newcommand{\bfg}{\boldsymbol{g}}
\newcommand{\ner}{\bol{r}}
\newcommand{\N}{\mathbb{N}}
\newcommand{\R}{\mathbb{R}}

\newcommand{\Rcal}{\mathcal R}
\newcommand{\Bcal}{\mathcal B}
\newcommand{\Ccal}{\mathcal C}
\newcommand{\Dcal}{\mathcal D}
\newcommand{\Pcal}{\mathcal P}
\newcommand{\Scal}{\mathcal S}
\newcommand{\Hcal}{\mathcal H}

\DeclarePairedDelimiter{\ceil}{\lceil}{\rceil}

\usepackage[bottom]{footmisc}

\counterwithout{theorem}{section}

\title{Construction of polynomial particular solutions of linear constant-coefficient partial differential equations}

\author{Thomas G. Anderson\thanks{Department of Computational Applied Mathematics \& Operations Research, Rice University, Houston, TX USA}, Marc Bonnet\thanks{POEMS (CNRS, INRIA, ENSTA), ENSTA Paris, 91120 Palaiseau, France}, Luiz M. Faria\textsuperscript{\textdagger}, and Carlos P\'erez-Arancibia\thanks{Department of Applied Mathematics, University of Twente, Enschede, The Netherlands.
\newline\phantom{aca}\hspace{.35em}Corresponding author: Marc Bonnet, e-mail: marc.bonnet@ensta-paris.fr}}

\date{\today}

\begin{document}

\maketitle

\begin{abstract}
    This paper introduces general methodologies for constructing closed-form solutions to linear constant-coefficient partial differential equations (PDEs) with polynomial right-hand sides in two and three spatial dimensions. Polynomial solutions have recently regained significance in the development of numerical techniques for evaluating volume integral operators and also have potential applications in certain kinds of Trefftz finite element methods. The equations covered in this work include the isotropic and anisotropic Poisson, Helmholtz, Stokes, linearized Navier-Stokes, stationary advection-diffusion, elastostatic equations, as well as the time-harmonic elastodynamic and Maxwell equations.  Several solutions to complex PDE systems are obtained by a potential representation and rely on the Helmholtz or Poisson solvers. Some of the cases addressed, namely Stokes flow, Maxwell's equations and linearized Navier-Stokes equations, naturally incorporate divergence constraints on the solution. This article provides a generic pattern whereby solutions are constructed by leveraging solutions of the lowest-order part of the partial differential operator (PDO). With the exception of anisotropic material tensors, no matrix inversion or linear system solution is required to compute the solutions.  This work is accompanied by a freely-available Julia library, \texttt{ElementaryPDESolutions.jl}, which implements the proposed methodology in an efficient and user-friendly format.
\end{abstract}

\section{Introduction}
This paper uses a combination of some simple ideas to obtain polynomial solutions of inhomogeneous constant-coefficient partial differential equations (PDEs) in $\mathbb R^d$ ($d=2,3$) with arbitrary given polynomial right-hand sides for many of the classical models arising in mathematical physics, in both the 2D and the 3D case. Such polynomial solutions hold in arbitrary regions, and are not constrained by conditions on a boundary or at infinity. The methods presented herein apply most basically to the familiar scalar partial differential operators (PDOs) but extend also to vector and anisotropic models. We detail the construction of polynomial solutions and publish an accompanying Julia library, \texttt{ElementaryPDESolutions.jl}\footnote{\href{https://github.com/IntegralEquations/ElementaryPDESolutions.jl/releases/tag/0.2}{https://github.com/IntegralEquations/ElementaryPDESolutions.jl}, version 0.2}.

While apparently simple and obviously not satisfactory for a complete theory or as a general method, polynomial solutions of the kind considered in this work have a demonstrated usefulness as components of other numerical solution techniques, such as the method of fundamental solutions and methods that use boundary integral equations for inhomogeneous PDEs. In the former, a particular solution is straightforwardly useful to reduce the problem to a homogeneous one that the method of fundamental solutions treats~\cite{Dangal2017}. In the latter approach, Green's identities are used to transform certain volume integrals to surface integrals, and in so doing particular solutions for PDEs corresponding to simple right-hand sides are introduced.  One of the first methods in this vein, dating to the 1980s, is the dual reciprocity method~\cite{partridge2012dual} that uses a global basis of simple functions such as monomials (alternatively, radial basis functions) to approximate an inhomogeneous right-hand side (sometimes called body force) in a linear PDE. The method thus calls for the associated (polynomial) solution to the PDE; see also~\cite{Atkinson1985} for another use of a global basis. Variants on the volume-to-boundary idea exist where the domain is meshed and approximation occurs on geometrically-simple regions~\cite{Greengard1996,Shen2022,Anderson2022}. Our interest arose in the course of using the latter kind of treatment as an indirect way to evaluate singular contributions to volume integral operators, and indeed the polynomial solutions given here allow the extension of the work~\cite{Anderson2022} to vectorial problems such as the Stokes, elasticity, and Maxwell systems.  Polynomial solutions also appear relevant to so-called Trefftz methods~\cite{cessenat1998application,Hiptmair2016} (and, possibly, related methods~\cite{LMIG2023} wherein polynomial solutions have been observed to form a basis with favorable conditioning properties); it will be of interest if the methods presented here are useful in such contexts.

Given a polynomial right-hand-side, the Helmholtz equation, and many other PDEs featuring a zeroth-order derivative term, have a unique polynomial solution. By contrast, such solution is defined up to arbitrary (e.g. harmonic) polynomials in the absence of a zeroth-order term, and in particular for the Poisson equation and other analogous cases such as elastostatics. Such solutions have been investigated for some time. Direct collocation approaches have long been used in the context of the method of fundamental solutions; see, e.g., \cite{poullikkas1998method,karageorghis2007efficient,gu2019localized}. Alternative solution methods which explicitly leverage the properties of polynomials have been developed to avoid the expense and ill-conditioning associated with certain linear systems of equations arising in collocation approaches. Recursions have been developed in~\cite{Janssen1992} for polynomial solutions to scalar constant-coefficient linear problems when $d=2$ or $d=3$; stability challenges in the recursion are discussed. Similar formulae for Poisson solutions are given in~\cite{Cheng1994} when $d = 2$, and this method is extended to the Helmholtz equation for $d=2$ and $d=3$ in~\cite{Golberg2003}. Solutions for polyharmonic and poly-Helmholtz operators are presented in~\cite{Tsai2009}. In all of these works, solutions are determined by utilizing a well-suited ansatz that relates the image of the relevant PDO to solutions corresponding to lower degree right-hand-sides.

In this work, we introduce and exploit general methodologies for constructing closed-form particular polynomial solutions to linear constant-coefficient PDEs for many classical models of mathematical physics. Our general goal is to obtain solutions for as many cases as possible by means of elementary and explicit algebraic methods only, in particular without having to (numerically) solve (often-underdetermined) linear systems verified by the coefficients of such polynomial solutions. The proposed methods, and resulting generality, hinge on the following aspects: (a) statement and systematic exploitation of the fact that solvability (in polynomials) of the lowest-order part of the relevant PDO always yields particular solutions of the whole PDE, (b) obtaining particular solutions of some commonly-involved lowest-order PDEs, and (c) exploiting the availability for more-complex (e.g. vector-valued) PDEs of representations by potentials to which solutions from (b) are applicable.

Perhaps the work closest to the present contribution can be found in~\cite{Dangal2017} for constant-coefficient second-order operators with a zeroth-order term present, as both express the solution in the form of the formal Neumann series expansion  $(k^2 + \Delta)^{-1}f=\sum_{j=0}^{\infty}(-1)^j k^{-2(j+1)}\Delta^j f$, $k\neq0$, where $f$ denotes the polynomial right-hand-side. Such an expansion gives rise to a finite number of terms by virtue of the fact that $\Delta$ is nilpotent as an operator on polynomials. Following these ideas, our Poisson solution approach relies on expressing the right-hand-side $f$ in terms of homogeneous polynomials and seeking $\Delta^{-1}f$ in the form $\sum_{j=0}^\infty c_j |\ner|^{2(j+1)}\Delta^jf$ with $\ner$ denoting the position vector, from which a simple recursion relation for the finite number of non-zero coefficients $c_j$ can be derived by applying Euler's theorem for homogeneous functions. Our approach to the Laplace operator appears to be novel and carries some advantages; unlike certain recurrence-based methods for obtaining a (non-unique) Poisson solution, the method described yields a solution that is much more symmetric in the input variables. Note that, unlike the unique polynomial solution of the Helmholtz equation, the Poisson polynomial solutions obtained by means of different methods need not necessarily coincide. That being said, we discovered after writing this article that the Poisson solution derived here was previously obtained in~\cite{Karachik_2010}, following an approach based on the analytical evaluation of volume potentials. We find the differing methodology employed herein to be readily generalizable, in particular to anisotropic models (i.e., PDOs involving $\operatorname{div}(A\nabla)$ where $A\in\mathbb C^{d\times d}$); of course, such solutions are naturally asymmetric.

Going beyond the Laplace and Helmholtz cases, the present paper constructs polynomial solutions to \emph{any} constant-coefficient general PDO assuming that solution techniques are available for the lowest-order part of the operator. A similar idea was used in~\cite{Dangal2017} in the specific case of the Helmholtz equation (wherein the lowest-order part of the operator is a constant times the identity and is easy to ``invert''), but the technique was not immediately generalizable and extension to other PDEs was considered an open research topic. The general expression for the polynomial solution presented here is inspired by the formal Neumann series expansion discussed above, again using nilpotence properties to result in a finite sum. Our method is demonstrated to obtain solutions for (anisotropic) advection-diffusion operators (in which the lowest order PDO is of first degree) as well as, in service of a real problem of interest arising in fluids, to a $6$\textsuperscript{th}-degree operator with a Laplacian as the lowest-order operator.

Known solution techniques for general classes of vectorial problems appear to be much more limited; to the best of the authors' knowledge, only the extension~\cite{Matthys1996} of the recursive technique~\cite{Janssen1992} exists for the solution of a certain class of vectorial problems. Although aimed at tackling quite general coupled vectorial PDE systems, such as the elastostatics system addressed in detail therein, the method presented in~\cite{Matthys1996} does not account for some important systems that feature divergence constraints on the PDE solution. In the case of the equations for two-dimensional Stokes flow, for instance, the PDE system can be written as
\begin{align}
\label{eq:stokes-explicit}
\begin{bmatrix}
\mu & 0 & 0 \\
0 & \mu & 0 \\
0 &0 & 0
\end{bmatrix}
\begin{bmatrix}
u_{xx}\\
v_{xx}\\
p_{xx}
\end{bmatrix}
+
\begin{bmatrix}
\mu & 0 & 0 \\
0 & \mu & 0 \\
0 &0 & 0
\end{bmatrix}
\begin{bmatrix}
u_{yy}\\
v_{yy}\\
p_{yy}
\end{bmatrix}
+
\begin{bmatrix}
0 & 0 & -1 \\
0 & 0 & 0 \\
1 & 0 & 0
\end{bmatrix}
\begin{bmatrix}
u_{x}\\
v_{x}\\
p_{x}
\end{bmatrix}
+
\begin{bmatrix}
0 & 0 & 0 \\
0 & 0 & -1 \\
0 & 1 & 0
\end{bmatrix}
\begin{bmatrix}
u_{y}\\
v_{y}\\
p_{y}
\end{bmatrix}
=
\begin{bmatrix}
f_1\\
f_2\\
0
\end{bmatrix},
\end{align}%
in the format used in~\cite{Matthys1996}, where $u,v$ denote the components of the velocity field, $p$ is the pressure, and $\mu$ is the viscosity. Since none of the four matrices in~\eqref{eq:stokes-explicit} are invertible, the determinant condition in~\cite[Sec.\ 3.2]{Matthys1996} is not satisfied and thus the recursive approach developed therein is not directly applicable to this equation.
Here, by contrast, we resort to our Helmholtz or Poisson solutions together with suitable representations of solutions of vectorial PDEs by potentials, to produce polynomial solutions that satisfy the unconstrained elastodynamic or elastostatic PDE systems as well as the constrained Maxwell or Stokes systems. The polynomial solutions we provide for the latter cases represent to our knowledge the first polynomial solutions to vectorial PDEs that fulfill the physically correct constraints on the vector field divergence, namely, the charge conservation and the incompressibility conditions in the cases of Stokes flow and Maxwell's equations, respectively. We additionally provide solutions to the Brinkman (a.k.a. linearized Navier-Stokes) system that also features an incompressibility constraint. Finally, we show that the polynomial solutions obtained for the anisotropic Poisson equation can be used to obtain solutions for some cases of anisotropic elastostatics.

Throughout this paper, we let $\Pcal_N$, $N\in\mathbb N_0$, denote the space of all $d$-variate polynomials of total degree at most $N$, with $d=2,3, \ldots$ the ambient dimension (a noteworthy fact being that with the exception of Maxwell all methods below, as well as the implementation, are dimension-agnostic). We make use of the standard multi-index notation where, for any $\alpha = (\alpha_1,\ldots,\alpha_d)\in\mathbb N_0^d$, we set $|\alpha| = \alpha_1+\ldots+\alpha_d$ and $\ner^\alpha=r_1^{\alpha_1} \ldots r_d^{\alpha_d}$ when $\ner=(r_1,\ldots,r_d)$. One simple but key fact that will be used in the sequel is that for any non-zero polynomial $p\in\Pcal_N$ and any PDO $\Bcal$ lacking a zeroth-order term, there exists a finite integer 
$m=m(p,\Bcal)$
such that $\Bcal^m p\not=0$ and $\Bcal^{m+1} p=0$.

The rest of this paper is organized as follows. In Section~\ref{sec:gen}, we present a general method to tackle linear constant-coefficient PDEs with non-homogeneous symbol by leveraging the assumed availability of a solution operator for the lowest-order term. We then proceed in Sections~\ref{sec:polynomial_solns:H} to~\ref{sec:polynomial_solns:L} to obtain particular polynomial solutions for PDEs with various lowest-order terms, for which we provide that solution operator; this in turn allows us to present additional solutions, pertaining to vectorial PDEs, which are enabled via suitable use of representations by potentials that satisfy previously-treated scalar PDEs. Finally, the accompanying Julia library and related numerical considerations are presented in Sec.~\ref{sec:numerical_aspects}.

\section{Preliminaries}\label{sec:gen}

In this paper, we investigate the problem of constructing polynomial solutions $u$ of inhomogeneous linear partial differential equations of the form
\begin{equation}
  \Bcal u = f \quad \text{in}\quad \R^d, \label{generic:pde}
\end{equation}
where the right-hand side $f$ is a given polynomial and the partial differential operator (PDO) $\Bcal$ has constant coefficients.
Any such PDO has the form $\Bcal=B(\partial)$, where $B$ is the total symbol of $\Bcal$, i.e., a $d$-variate polynomial of degree $p$, with $p$ thus being the order of $\Bcal$. More precisely:

\begin{definition}
Let $\Hcal_n\subset\Pcal_n, n \in \N_0,$ denote the $d$-variate homogeneous polynomials of total degree~$n$.
By the \emph{total symbol} of the constant-coefficient partial differential operator $\Bcal:=B(\partial)$ of order $p\in\mathbb N$, with $p\geq q\in\mathbb{N}_0$, we shall mean
$$
B(\bol\xi)=\sum_{j=q}^{p} a_jB_j(\bol\xi),\quad\bol\xi\in\R^d,
$$
where $B_j\in \Hcal_j$
and $a_j\in\mathbb C$, with $a_q \ne 0$ and $a_p \ne 0$.
For any $j \in \{q,\ldots,p\}$, $B_j$ will be referred to as the \emph{homogeneous part of degree $j$} of $B$; $B_p$ is called the \emph{principal (or leading) part} of the symbol $B$ and $B_q$ is called the \emph{elementary part} of the symbol $B$ (or the \emph{elementary symbol}).
\end{definition}

Clearly, whenever the symbol $B$ is not itself homogeneous, we can write $B=a_qB_{q}+R_q$, where $a_q B_{q}\in\Hcal_q$ is the elementary part of $B$ (having degree $q$) and $  R_q:=\sum_{j=q+1}^p a_j B_j\in\Pcal_p$ defines the remainder (higher-degree) symbol, with degree $p>q\geq0$.  
We additionally denote by $r$ the first integer greater than $q$ for which $a_r \ne 0$; $r$ is thus the degree of the elementary symbol of $ \Rcal_q$ and $r-q$ measures the minimum difference in derivative orders extant in $\Bcal_q$ and $ \Rcal_q$. With these definitions, the PDO $\Bcal$ in the generic problem~\eqref{generic:pde} has the additive decomposition
\begin{equation}
  \Bcal = \Bcal_q + \Rcal_q  \label{B:split}
\end{equation}
where the PDOs $\Bcal_q = a_q B_q(\partial)$ and $ \Rcal_q = R_q(\partial)$ correspond to the elementary and remainder symbols.

Our methodology rests on the assumption that a polynomial solution $v\in\Pcal_{n+q}$ of 
\begin{equation}
  \Bcal_q v = g \quad \text{in}\ \ \R^d,\label{generic:split2}
\end{equation}
can be constructed for any given $g\in\Pcal_n,\,n\in\{0,\ldots,N\}$ and so define $\Scal_q$ as a \emph{right inverse} of $\Bcal_q$ such that $\Scal_q$ is linear, $v=\Scal_q g \in \Pcal_{n+q}$. In particular, we provide this right inverse explicitly for certain $\Bcal_q$, $q = 0, 1, 2, 4$, in what follows. Our PDE polynomial solution construction methodology for the original problem~\eqref{generic:pde} then makes use of the following result:
\begin{theorem}\label{thm:soln_expr}
Let $\Bcal$ be decomposed as in~\eqref{B:split}, and let $\Scal_q$ denote a right inverse of $\Bcal_q$ having the properties specified above.  For any given non-zero $f\in\Pcal_N$, the polynomial $u\in\Pcal_{N+q}$ given by 
\begin{equation}
u = \Scal_q \sum_{j=0}^{m} \big( -{\Rcal}_q\Scal_q \big)^jf,
\label{eq:neumannseries_generic_repr}
\end{equation}
with $m\in\N_0$ such that $\big( {\Rcal}_q\Scal_q \big)^m f\not=0$ and $\big( {\Rcal}_q\Scal_q \big)^{m+1}f=0$, satisfies
\begin{equation}
  \Bcal u = f \quad \text{in}\ \ \R^d \label{generic:split}.
\end{equation}
Thus $\Scal:=\Scal_q \sum_{j\geq 0} \big( -{\Rcal}_q\Scal_q \big)^j$ defines a right inverse of $\Bcal$ on polynomials, the sum being always finite (cf.\ \eqref{eq:neumannseries_generic_repr}); in particular\footnote{Because $\Scal_q 0 = 0$, the precise value $m$ need not be identified; in fact, for simplicity the implementation computes $u$ via the sum in~\eqref{eq:neumannseries_generic_repr} up to the provided upper bound for $m$.}  $m \le \ceil*{(N+1)/(r - q)}-1$.

In addition, we have
\begin{equation}
  \text{(a) \ }\operatorname{ker}(\Bcal|_{\Pcal_N})=\big\{ w-\Scal\Rcal_q w,\, w\in\operatorname{ker}(\Bcal_q|_{\Pcal_N}) \big\} \qquad
  \text{and} \quad \text{(b) \ }\operatorname{dim}\operatorname{ker}(\Bcal|_{\Pcal_N})=\operatorname{dim}\operatorname{ker}(\Bcal_q|_{\Pcal_N}), \label{nullspace}
\end{equation}
where $\Dcal|_{\Pcal_N}$ denotes a generic linear, constant-coefficient PDO $\Dcal$ restricted to $\Pcal_N$. 
In particular, if $q=0$ and the
multiplication operator $\Bcal_0$ is invertible, we have $\operatorname{ker}(\Bcal|_{\Pcal_N})=\operatorname{ker}(\Bcal_0|_{\Pcal_N})=\{0\}$, so that~\eqref{eq:neumannseries_generic_repr} defines the unique polynomial solution of~\eqref{generic:pde}.
\end{theorem}
\begin{proof}
Using the splitting $\Bcal=\Bcal_q+\Rcal_q$ we have 
$$
  \Bcal\Scal_q \sum_{j=0}^{m} \big( -{\Rcal}_q\Scal_q \big)^jf
 = \sum_{j=0}^{m} \left\{ \big( -{\Rcal}_q\Scal_q \big)^jf - \big( -{\Rcal}_q\Scal_q \big)^{j+1}f \right\}
 = f + (-1)^m \big( {\Rcal}_q\Scal_q \big)^{m+1}f = f,
$$
where the next-to-last equality above is obtained from the telescoping sum and the last results from the definition of $m$. To prove the claimed (finite) upper bound on $m$, we note that 
$\Bcal^{-1}_q:\Pcal_{n}\to\Pcal_{n+q}$ for any $n\in\mathbb N_0$, and hence ${\Rcal}_q\Scal_q:\Pcal_{n}\to\Pcal_{n+q-r}$, implying that $\big( {\Rcal}_q\Scal_q \big)^{\ell+1}f\in\Pcal_{N+(\ell+1)(q-r)}$ for $\ell \in \N_0$. With $\ell=\ceil*{(N+1)/(r - q)}-1$, we thus have $N+(\ell+1)(q-r)<0$, hence $\big( {\Rcal}_q\Scal_q \big)^{\ell+1}f=0$ resulting in $m\leq\ell$. Finally, the polynomial $u$ so defined satisfies $u \in \Pcal_{N+q}$ by the aforementioned mapping properties of $\Rcal_q$ and $\Scal_q$, the highest possible degree contribution to $u$ occurring in the sum in~\eqref{eq:neumannseries_generic_repr} at $j = 0$.

To prove the claims in~\eqref{nullspace}, we begin by 
decomposing  $u\in\Pcal_N$ as $u = u_N + v$ with $u_N\in\Hcal_N$ and $v\in\Pcal_{N-1}$. Since $\Bcal_q u_N\in\Hcal_{N-q}$ while $(\Bcal_q v + \Rcal_q u)\in\Pcal_{N-q-1}$, a necessary condition for $u$ to satisfy $\Bcal u = (\Bcal_q + \Rcal_q)u = 0$ is $u_N\in \operatorname{ker}(\Bcal_q|_{\Hcal_N})$.
Then, $u=u_N+v\in\operatorname{ker}(\Bcal|_{\Pcal_N})$ if and only if $\Bcal v=-\Rcal_q u_N$. The solutions in $\Pcal_{N-1}$ of the latter equation are all $v$ of the form $v=-\Scal\Rcal_q u_N+w$ with $w\in\operatorname{ker}(\Bcal|_{\Pcal_{N-1}})$, where $\Scal$ is the right inverse of $\Bcal$ defined by~\eqref{generic:split}. Any $u=u_N+v\in\operatorname{ker}(\Bcal|_{\Pcal_N})$ is therefore found to have the form $u=(I-\Scal\Rcal_q) u_N + w$ with $u_N\in \operatorname{ker}(\Bcal_q|_{\Hcal_N})$ and $w\in\operatorname{ker}(\Bcal|_{\Pcal_{N-1}})$.

The foregoing argument reduces the characterization of $\operatorname{ker}(\Bcal|_{\Pcal_N})$ to that of $\operatorname{ker}(\Bcal|_{\Pcal_{N-1}})$, and can hence be applied recursively to $\operatorname{ker}(\Bcal|_{\Pcal_\ell})$ along decreasing degrees $\ell=N,N-1,\ldots,1$. This results in any $u\in\operatorname{ker}(\Bcal|_{\Pcal_N})$ having the representation
\begin{equation}
  u = (I-\Scal\Rcal_q)\sum_{\ell=1}^N u_{\ell} \qquad \text{with \ }  u_{\ell}\in \operatorname{ker}(\Bcal_q|_{\Hcal_{\ell}}).
\end{equation}
In addition, we have that $\operatorname{ker}(\Bcal_q|_{\Pcal_N})=\bigoplus_{\ell=0}^{N}\operatorname{ker}(\Bcal_q|_{\Hcal_{\ell}})$  (by virtue of $B_q$ being homogeneous). The above representation of $u\in\operatorname{ker}(\Bcal|_{\Pcal_N})$ therefore yields the characterization~(\ref{nullspace}a). 

Finally, we note that $(I-\Scal\Rcal_q) u_{\ell}=0$ with $u_{\ell}\in \operatorname{ker}(\Bcal_q|_{\Hcal_{\ell}})$ requires $u_{\ell}=0$. This remark, together with the above direct-sum form of $\operatorname{ker}(\Bcal_q|_{\Pcal_N})$, implies that the null space of $I-\Scal\Rcal_q$ acting on $\operatorname{ker}(\Bcal_q|_{\Pcal_N})$ is trivial, and the claim~(\ref{nullspace}b) follows.
\end{proof}

Theorem~\ref{thm:soln_expr} emphasizes the essential role played by the elementary symbol $B_q$ for finding polynomial solutions to~\eqref{generic:pde}.
In the remainder of this paper, we address various classical PDEs, proceeding along increasing degree of 
$B_q$ (increasing values of the lowest order $q$ present in $\Bcal$). 
For $q=0$ (e.g.\ where the elementary part of $\Bcal$ is a multiplicative scalar), the archetypal situation is the scalar Helmholtz equation (where $\Bcal=\Delta+k^2$, $\Bcal_0=k^2$, $\Rcal_0=\Delta$ and $p = r = 2$); this case is addressed in Sec.~\ref{sec:polynomial_solns:H}. Then, we solve in Sec.~\ref{sec:diff:advec} the case where the elementary part $\Bcal_1$ is the advection operator $\bol{\beta}\cdot\nabla$, for which $q=1$. Finally, we consider in Sec.~\ref{sec:polynomial_solns:L} cases where the elementary part $\Bcal_q$ is the (isotropic or anisotropic) Laplacian ($q=2$) or the isotropic bilaplacian ($q = 4$).

\section{PDOs with zeroth-order elementary part (\emph{\`a la} Helmholtz)}\label{sec:polynomial_solns:H}

\subsection{Helmholtz\label{sec:scalar_Helm}}

We start off by considering  the Helmholtz equation with wavenumber $k\in\mathbb C$, $k\neq0$.
The problem at hand is thus: given $f\in \Pcal_N$, find $u\in\Pcal_N$ such that
\begin{equation}
  (\Delta+k^2)u = f\quad\text{in}\quad \R^d.  \label{helmholtz:pb}
\end{equation}

Here, the lowest-order generic problem~\eqref{generic:split2} is the simple zeroth-order equation $k^2 v = g$, whose unique solution is of course $v = k^{-2} g$. The solution of~\eqref{helmholtz:pb} is hence unique (by Theorem~\ref{thm:soln_expr}), and is provided by~\eqref{eq:neumannseries_generic_repr} as 

\begin{equation}
  u = \sum_{j=0}^m (-1)^jk^{-2(j+1)}\Delta^j f
  = k^{-2}f - k^{-4}\Delta f + k^{-6}\Delta^2 f - \ldots -(-k^2)^{-m-1}\Delta^m f. \label{P:sol:helm}
\end{equation}

\subsection{Elastodynamics\label{sec:somigliana}} Let us now consider the (linear, isotropic) time-harmonic elastodynamic system of equations corresponding to a medium endowed with mass density $\rho>0$, shear modulus $\mu>0$, and Poisson's ratio $\nu\in(0,\frac12)$.
The problem at hand in this case is: given $\bol f\in [\Pcal_N]^d$, find $\bol u\in[\Pcal_N]^d$ satisfying the Navier equation
\begin{equation}
   \Delta\neu  + \frac{1}{(1-2\nu)}\nabla(\operatorname{div}\neu)+  k_2^2\neu  = \nef\quad\text{in}\quad\R^d,    \label{eq:elastodyn:vector}
\end{equation}
where  $k_2=\omega/c_2$ is the shear wavenumber defined in terms of the velocity $c_2=\sqrt{\mu/\rho}$ and angular frequency $\omega>0$.

It is well known that any time-harmonic displacement field $\neu$ expressed in the form~\cite[Sec.~5.4]{eri}
\begin{equation}
  \neu = 2(1-\nu)(\Delta+k^2_1)\bfg - \nabla(\operatorname{div}\bfg) \label{u:Somi}
\end{equation}
satisfies the Navier equation~\eqref{eq:elastodyn:vector} with given body force density $\mu\nef$ provided the \emph{Somigliana vector potential} $\bfg$ satisfies the repeated vector Helmholtz equation
\begin{equation}
  (\Delta+k^2_1)(\Delta+k^2_2)\bfg  = \frac{\nef}{2(1-\nu)}\quad\text{in}\quad\R^d \label{eq:elastodyn}
\end{equation}
(since substituting the ansatz~\eqref{u:Somi} into~\eqref{eq:elastodyn:vector} yields~\eqref{eq:elastodyn}), where $k_1=\omega/c_1$ is the compressional wavenumber, the corresponding velocity being given by $c_1=\sqrt{2\mu(1-\nu)/\rho(1-2\nu)}$.

Given a vector-valued polynomial $\nef\in [\Pcal_N]^d$, we then obtain the polynomial solution $\bfg$ of~\eqref{eq:elastodyn} by consecutively applying the vector-valued version of~\eqref{P:sol:helm} to the following  two  inhomogeneous componentwise-scalar Helmholtz problems:
\begin{equation}
\begin{aligned}
  &\text{(a) \ find } \bol{q}\in[\Pcal_N]^d \quad\text{such that}\qquad &(\Delta+k^2_1)\bol{q} &=\frac{\nef}{2(1-\nu)}, \\
  &\text{(b) \ find } \bfg\in[\Pcal_N]^d \quad\text{such that} &(\Delta+k^2_2)\bfg &= \bol{q}.
\end{aligned} \label{eq:somi:hel}
\end{equation}
The displacement solution $\neu\in[\Pcal_N]^d$, which by Theorem~\ref{thm:soln_expr} (see also Sec.~\ref{PDO:lowerorder}) is the unique polynomial solution to~\eqref{eq:elastodyn:vector}, is then found to be explicitly given by~\eqref{u:Somi} with the above solution $\bfg\in [\Pcal_N]^d$, which is obtainable for any given $\bol f\in [\Pcal_N]^d$.

The foregoing procedure highlights the potential usefulness of available representations by potentials of solutions of e.g. vector PDEs.
Equation~\eqref{eq:elastodyn:vector} can alternatively be solved by means of a direct application of Theorem~\ref{thm:soln_expr} to either the original PDE~\eqref{eq:elastodyn:vector} (with $\bol{\Bcal}_0 = k_2^2\bol{I}$ and $\bol{\Rcal}_0 = \Delta + (1-2\nu)^{-1}\nabla(\operatorname{div}$) or to the equations~\eqref{eq:elastodyn} governing the vector potential (with $\Bcal_0 = k_1^2k_2^2$ and $\Rcal_0 = \Delta^2 + (k_1^2 + k_2^2)\Delta$).

\subsection{General linear PDO having a zeroth-order elementary part}
\label{PDO:lowerorder}

For a given scalar $c\not=0$, the generic equation
\begin{equation}
  \Rcal_0 u + c u =  f\quad\text{in}\quad \R^d,
\end{equation}
(corresponding to $\Bcal_0=c$ with the notations of Sec.~\ref{sec:gen})

has for a polynomial right-hand side $ f$ a unique
polynomial solution $u$ given by
\begin{equation}
  u = c^{-1}f - c^{-2}\Rcal_0 f + c^{-3}(\Rcal_0)^2 f - c^{-4}(\Rcal_0)^3 f \ldots, \label{sol:order0}
\end{equation}
where the sum

terminates after a finite number of terms.
This specialization of Theorem~\ref{thm:soln_expr} coincides in the case $p = 2$ with the method given in~\cite{Dangal2017}
.
The solution expression~\eqref{sol:order0} extends straightforwardly to vector-valued PDOs of the form $\bol{\Rcal}_0\neu + \bol{\Bcal}_0\neu$, where $\bol{\Rcal}_0$ is a matrix-valued PDO with no zeroth-order term and 
the elementary part $\bol{\Bcal}_0\in\mathbb{C}^{d\times d}$ is a matrix-valued constant. This includes e.g. the elastodynamics operator for general anisotropic elastic media, and in particular provides an alternative to the method of Sec.~\ref{sec:somigliana} for isotropic elastodynamics (wherein $\bol{\Bcal}_0=k_2^2 \bol{I}$).

We observe that the solution~\eqref{sol:order0} for $c\not=0$ cannot be used to obtain a solution of $\Rcal_0 u = f$ (e.g.\ $\Rcal_0 = \Delta$) by passing to the limit $c\to 0$ as the former solution is of degree $N$ whereas any solution of the latter problem must be of degree at least $N+r$, so cannot be attained by such putative limiting processes. In addition, the solution~\eqref{sol:order0} for fixed $f$ clearly blows up as $c\to 0$.\bigskip

\subsection{Maxwell}\label{sec:maxwell}

The corresponding polynomial problem for the Maxwell equation system consists in finding the time-harmonic electromagnetic field $(\bol{E},\bol H)\in [\Pcal_N]^3\times[\Pcal_{N-1}]^3$, arising in a homogeneous isotropic medium (with constant scalar permittivity $\varepsilon$ and permeability $\mu$) and due to a given polynomial current $\bol{J}\in [\Pcal_N]^3$ and an associated polynomial charge density $\rho\in\Pcal_{N-1}$, which satisfies:
\begin{equation}
\begin{aligned}
  \mathrm{i}\omega\varepsilon\bol{E} + \operatorname{rot}\bol{H} &= \bol{J}, \qquad &
  -\mathrm{i}\omega\mu\bol{H} + \operatorname{rot}\bol{E} &= \bol{0}, \\
  \varepsilon\operatorname{div}\bol{E} &= \rho, &
  \mu\operatorname{div}\bol{H} &= 0,
\end{aligned} \label{maxwell}
\end{equation}
in $\R^3$, with the sources being constrained by the charge conservation equation:
\begin{equation}
  \operatorname{div}\bol{J} - \mathrm{i}\omega\rho = 0\quad\text{in}\quad\R^3. \label{charge:conserv}
\end{equation}

We seek the electromagnetic fields $(\bol{E},\bol{H})$ in the well-known potential-representation form~\cite{Stratton:41}
\begin{equation}
  \bol{E} = \mathrm{i}\omega\bol{A} - \nabla\varphi, \qquad \bol{H} = \frac{1}{\mu}\operatorname{rot}\bol{A}, \label{EH:pot}
\end{equation}
where the vector potential $\bol{A}$ and the scalar potential~$\varphi$ are in addition constrained through a gauge condition to prevent indetermination. Here we assume $(\bol{A},\varphi)$ to be linked by the Lorenz gauge condition:
\begin{equation}
  \operatorname{div}\bol{A} - \mathrm{i}\omega\varepsilon\mu\varphi = 0\quad\text{in}\quad\R^3, \label{lorenz}
\end{equation}
Upon inserting the ansatz~\eqref{EH:pot} in~\eqref{maxwell} and taking the gauge condition~\eqref{lorenz} into account, the fields $(\bol{E},\bol{H})$ so represented satisfy Maxwell's equations~\eqref{maxwell} provided $\bol{A}$ and $\varphi$ solve the Helmholtz problems
\begin{equation}
\begin{aligned}
& \text{(a) find }\bol A\in [\Pcal_N]^3 \quad\text{such that}  \qquad &(\Delta+k^2)\bol{A} &= - \mu\bol{J} , \\
& \text{(b) find }\varphi\in\Pcal_{N-1} \,\quad\text{such that}  &(\Delta+k^2)\varphi &=- \rho/\varepsilon , \label{wave:potentials}
\end{aligned}
\end{equation}
where $k=\omega\sqrt{\varepsilon\mu}$ is the wavenumber. As shown in Sec.~\ref{sec:scalar_Helm}, the unique pair of polynomial potentials solving the Helmholtz equations~\eqref{wave:potentials} with arbitrary polynomial source terms are given by
\begin{equation}
  \bol{A} =-\mu\sum_{j=0}^m(-k^2)^{-j-1}\Delta^j\bol J\ \in[\Pcal_N]^3,\qquad
  \varphi =-1/\varepsilon\sum_{j=0}^m(-k^2)^{-j-1}\Delta^j \rho\ \in\Pcal_{N-1}.\label{eq:maxwell_soln_rep}
\end{equation}
Moreover, it is easy to show that the polynomial solution $(\bol A,\varphi)$ given by~\eqref{eq:maxwell_soln_rep} does satisfy the Lorenz gauge condition~\eqref{lorenz}. Indeed, substituting the solutions~\eqref{eq:maxwell_soln_rep} in~\eqref{lorenz} and in view of the charge conservation equation~\eqref{charge:conserv},
we have
$$
  \operatorname{div}\bol{A} - \mathrm{i}\omega\varepsilon\mu\varphi
 = -\mu\sum_{j=0}^m(-k^2)^{-j-1}\Delta^j\left(\operatorname{div}\bol J - \mathrm{i}\omega\rho \right)=0.
$$
In summary, the unique pair of polynomial potentials solving the wave equations~\eqref{wave:potentials} for arbitrary (polynomial) sources $(\bol A,\varphi)$ that verify the requisite charge conservation constraint~\eqref{charge:conserv} (where for instance $\bol{J}\in[\Pcal_N]^d$ is arbitrary and $\rho\in\Pcal_{N-1}$ is then given by~\eqref{charge:conserv}) also satisfies automatically the Lorenz gauge. The electromagnetic fields~\eqref{EH:pot} with $(\bol A,\varphi)$ given by~\eqref{eq:maxwell_soln_rep} are therefore the unique polynomial electromagnetic solution in $\R^3$ of the potential-representation form~\eqref{EH:pot}, \eqref{lorenz} created by given polynomial sources $\bol{J},\rho$.

In fact, those fields are the unique polynomial solution of the Maxwell system~\eqref{maxwell}. If not, any nontrivial polynomial solution of the homogeneous version of system~\eqref{maxwell} must (on combining the first two equations in two different ways) satisfy $\operatorname{rot}\operatorname{rot}\bol{H}-k^2\bol{H}=\bol{0}$ and $\operatorname{rot}\operatorname{rot}\bol{E}-k^2\bol{E}=\bol{0}$. The zeroth-order elementary part of the latter equations being invertible, they (and hence the homogeneous system~\eqref{maxwell}) admit only the trivial polynomial solution (see also Sec.~\ref{PDO:lowerorder}), from which the claimed uniqueness follows.

\section{PDOs with first-order elementary part (\emph{\`a la} advection)}\label{sec:diff:advec}

We next seek polynomial solutions $u$, for given $f\in\Pcal_{N}$, of equations of the form
\begin{equation}
  \Rcal_1 u + \bol{\beta} \cdot \nabla u = f, \label{diff:advec}
\end{equation}
where $q=1$, $\Rcal_1$ only contains derivatives of at least second order (in particular, $p \ge r \ge 2$)
and $\Bcal_1:=\bol{\beta} \cdot \nabla$ is the advection operator with $\bol{\beta}\in\mathbb{R}^d$ an arbitrary non-zero constant vector. For instance, the diffusion-advection equation has the form~\eqref{diff:advec} with $\Rcal_1$ taken as the (isotropic or anisotropic) Laplacian.

The procedure described in Sec.~\ref{sec:gen} applies provided a particular solution $v\in\Pcal_{n+1}$ of
\begin{equation}
  \bol{\beta} \cdot \nabla v = g \label{advec_only}
\end{equation}
can be found for any given $g\in\Pcal_n$, $n=0,\ldots,N$. To this end, $v$ is sought of the form
\begin{equation}
  v = \frac{1}{|\bol{\beta}|^2} \sum_{\ell=0}^n c_\ell\, (\bol{\beta}\cdot\bol{r})^{\ell+1}\, \Bcal_1^{\ell}g.
\label{diff:advec:ansatz}
\end{equation}
Evaluating $\Bcal_1 v$ on the above ansatz and rearranging, we find
\begin{equation}
  \Bcal_1 v = \bol{\beta} \cdot \nabla v = \frac{1}{|\bol{\beta}|^2}  \Big\{ c_0 |\bol{\beta}|^2 g
  + \sum_{\ell=0}^{n-1} \big[\, c_{\ell+1} (\ell+2) |\bol{\beta}|^2 + c_{\ell} \big] \, \Bcal_1^{\ell+1}g \Big\},
\end{equation}
so that setting
\begin{equation}
  c_0 = 1, \qquad\quad c_{\ell+1} = -\frac{c_{\ell}}{(\ell+2) |\bol{\beta}|^2} \quad (0\leq \ell\leq n-1), \label{diff:advec:ansatz:csts}
\end{equation}
in~\eqref{diff:advec:ansatz} produces a polynomial $v\in\Pcal_{m+1}$ that solves~\eqref{advec_only}. Consequently, a particular polynomial solution of~\eqref{diff:advec} is obtained from formula~\eqref{eq:neumannseries_generic_repr} with $\Scal_1$ given by~\eqref{diff:advec:ansatz} and the constants $c_{\ell}$ prescribed in~\eqref{diff:advec:ansatz:csts}.

\section{PDOs with higher-order elementary part (\emph{\`a la} Laplace)}\label{sec:polynomial_solns:L}

\subsection{Poisson/Laplace\label{sec:Laplace}} In this section we focus on the polynomial Poisson problem: given $f\in \Pcal_N$, find a polynomial $u\in\Pcal_{N+2}$ such that
\begin{equation}
  \Delta u = f \quad\text{in}\quad \R^d. \label{poisson:pb}
\end{equation}
In terms of the definitions of Section~\ref{sec:gen}, we have $\Bcal=\Bcal_q=\Delta$ ($\Bcal$ being homogeneous) and $\Rcal_q=0$.
We propose a generic solution method for~\eqref{poisson:pb} based on the following observations:\enlargethispage*{10ex}
\begin{itemize}
\item Any $f\in \Pcal_N$ is a finite sum of homogeneous polynomials of degree at most $N$. By linearity, we can then simplify the problem~\eqref{poisson:pb} by assuming that $f\in\Hcal_n$. Note that $\operatorname{dim}(\Hcal_n)=n+1$ if $d=2$ and $\operatorname{dim}(\Hcal_n)=(n+1)(n+2)/2$ if $d=3$.
For example, the Maclaurin series
$$
f(\ner)= \sum_{|\alpha|=0}^N \frac{D^\alpha f(\bol 0)}{\alpha!} \ner^\alpha,\quad \ner\in\R^d,
$$
provides the expansion of a given polynomial $f\in\Pcal_N$ on the basis of homogeneous polynomials $\ner^\alpha\in \Hcal_{|\alpha|}$ $(0\leq|\alpha|\leq N)$.

\item For some $\ell\in\mathbb N_0$ and $n\in\mathbb N_0$, let $v=r^{2\ell+2}h_n$ with $h_n\in\Hcal_n$ and $r=|\ner|$. Then:
\begin{equation}
  \Delta v = \gamma_{\ell}^{n}r^{2\ell}h_n + r^{2\ell+2} \Delta h_n\quad\text{and}\quad \gamma_{\ell}^{n} := 2(\ell+1)(2\ell+2n+d). \label{DeltaP}
\end{equation}
Indeed, elementary computations yield
\begin{equation}
  \nabla(r^{2\ell+2}) = (2\ell+2) r^{2\ell}\ner\quad\text{and}\quad
  \Delta (r^{2\ell+2}) = 2(\ell+1)(2\ell+d)\,r^{2\ell}. \label{Delta:ids}
\end{equation}
Therefore,
\[
\begin{aligned}
  \Delta v &= (\Delta r^{2\ell+2})h_n+r^{2\ell+2}(\Delta h_n)+2\nabla(r^{2\ell+2})\cdot\nabla h_n \\
  &= 2(\ell+1)(2\ell+d)h_n+4(\ell+1)r^{2\ell}\ner\cdot\nabla h_n + r^{2\ell+2} \Delta h_n \\
  &=\gamma_{\ell}^{n}r^{2\ell}h_n + r^{2\ell+2} \Delta h_n,
\end{aligned}
\]
where we have used the fact that $\ner\cdot\nabla h_n=nh_n$, i.e., Euler's theorem for homogeneous functions.
\end{itemize}

Let then $f\in\Hcal_n$ (e.g. a monomial).  We seek a solution $u\in\Pcal_{n+2}$ of $\Delta u=f$ of the form
\[
  u = \sum_{\ell=0}^m c_\ell r^{2\ell+2} \Delta^\ell f,
\]
where the coefficients $\{c_\ell\}_{\ell=0}^m$ are to be determined and where $m=m(f)$ is defined, as before, as the smallest integer such that $\Delta^{m+1}f\equiv 0$. To this end, we first observe that $\Delta^\ell f\in \Hcal_{n-2\ell}$. Therefore, by virtue of~\eqref{DeltaP}, we have
\[
\begin{aligned}
  \Delta u
  &= \sum_{\ell=0}^m c_\ell \gamma_\ell^{n-2\ell} r^{2\ell} \Delta^\ell f + \sum_{\ell=0}^{m-1} c_\ell r^{2\ell+2} \Delta^{\ell+1} f \\
  &= c_0 \gamma_0^n f + \sum_{\ell=1}^m \big( \gamma_\ell^{n-2\ell} c_\ell + c_{\ell-1} \big)\, r^{2\ell}\Delta^\ell f.
\end{aligned}
\]
The polynomial $u$ thus satisfies $\Delta u=f$ if we set
\begin{equation}
  c_0 = \frac{1}{\gamma_0^n} = \frac{1}{2(2n+d)} \qquad\text{and}\qquad c_\ell = -\frac{c_{\ell-1}}{\gamma_\ell^{n-2\ell}} = -\frac{c_{\ell-1}}{2(\ell+1)(2n-2\ell+d)},\quad 1\leq\ell\leq m.
\label{Delta:ck}
\end{equation}

\subsection{Bilaplacian\label{sec:bilaplacian}}
We next consider the case of the bilaplacian PDO. The polynomial solutions developed in this section will become important in the sequel to construct solutions to the inhomogeneous elastostatic and Stokes equations. We seek a polynomial solution $u\in\Pcal_{N+4}$ of
\begin{equation}
  \Delta^2 u = f \quad\text{in}\quad \R^d, \label{polyn:bihar}
\end{equation}
where $f\in\Pcal_N$, and note that, following the construction of the Poisson equation solutions presented above in Sec.~\ref{sec:Laplace}, it suffices to restrict ourselves to the case of the source a homogeneous polynomial $f\in\Hcal_n$, $n\in\mathbb N_0$. Finding a particular $u\in \Hcal_{n+4}$ that solves~\eqref{polyn:bihar} given $f\in \Hcal_{n}$ is a straightforward task using either of the following two approaches.

The {\bf first approach} simply consists in applying twice the procedure for the Poisson equation, i.e.,
\[
  \text{(a) \ find $g\in\Hcal_{n+2}$ such that \ } \Delta g = f, \qquad
  \text{(b) \ find $u\in\Hcal_{n+4}$ such that \ } \Delta u = g.
\]

The {\bf second approach} consists in directly seeking $u\in\Hcal_{n+4}$ as the sum
\[
  u = \sum_{\ell=0}^m c_\ell r^{2\ell+4} \Delta^\ell f,
\]
so that it satisfies $\Delta^2 u = f$.  Then, in view of the identity
\[
  \Delta^2 r^{2\ell+4} h_n
 = \gamma^n_{\ell+1}\gamma^n_\ell\, r^{2\ell}h_n + 2\gamma^n_{\ell+1}r^{2\ell+2}\Delta h_n + r^{2\ell+4} \Delta^2 h_n,
\]
obtained by applying $\Delta$ to $\Delta r^{2\ell+4}h_n$ using twice the formula~\eqref{DeltaP} and exploiting the fact that $\Delta h_{n}\in\Hcal_{n-2}$ for $h_n\in\Hcal_n$, we arrive at
\[
\begin{aligned}
  \Delta^2 u
 &= \sum_{\ell=0}^m c_\ell \gamma_{\ell}^{n-2\ell} \gamma_{\ell+1}^{n-2\ell} r^{2\ell} \Delta^\ell f
  + \sum_{\ell=0}^{m-1} 2c_\ell \gamma_{\ell+1}^{n-2\ell} r^{2\ell+2} \Delta^{\ell+1} f
  + \sum_{\ell=0}^{m-2} c_\ell r^{2\ell+4} \Delta^{\ell+2} f \\
 &= c_0 \gamma_0^n \gamma_1^n f + \big( c_1 \gamma_1^{n-2} \gamma_2^{n-2} + 2c_0 \gamma_1^n \big) r^2\Delta f
 + \sum_{\ell=2}^m \big( c_\ell \gamma_{\ell}^{n-2\ell} \gamma_{\ell+1}^{n-2\ell} + 2c_{\ell-1}\gamma_\ell^{n+2-2\ell}  + c_{\ell-2} \big)  r^{2\ell} \Delta^\ell f.
\end{aligned}
\]
The polynomial $u\in\Hcal_{n+4}$ is thus made to satisfy $\Delta^2 u=f\in\Hcal_n$ by recursively defining the coefficients $\{c_\ell\}_{\ell=0}^m$ as
\[
  c_0 = \frac{1}{\gamma_0^n \gamma_1^n}, \quad  c_1 =-\frac{2\gamma_1^n c_0}{\gamma_1^{n-2} \gamma_2^{n-2}}, \qquad
  c_\ell= - \frac{2\gamma_\ell^{n+2-2\ell}c_{\ell-1}}{\gamma_{\ell}^{n-2\ell} \gamma_{\ell+1}^{n-2\ell}} - \frac{c_{\ell-2}}{\gamma_{\ell}^{n-2\ell} \gamma_{\ell+1}^{n-2\ell}}, \quad \ell\geq 2.
\]

The solutions $u\in\Hcal_{n+4}$ yielded by the two proposed approaches for solving $\Delta^2 u=f$ are not identical. For simplicity, the first approach for computing $u$ is implemented in \texttt{ElementaryPDESolutions.jl}.

\subsection{Isotropic elastostatics}
\label{elas:iso}
Let us now consider a (linearly elastic, isotropic) medium endowed with shear modulus $\mu>0$ and Poisson's ratio $\nu\in(0,\frac12)$. The problem at hand is then to find a elastostatic displacement field $\neu\in[\Pcal_{N+2}]^d$ generated by a given body force density $\mu\nef\in[\Pcal_N]^d$, which satisfies the PDE system
\begin{equation}
  \Delta\neu+ \frac{1}{(1-2\nu)}\nabla(\operatorname{div}\neu)   = \nef\quad\text{in}\quad\R^d. \label{PDE:elas}
\end{equation}

The zero-frequency form of the approach of Sec.~\ref{sec:somigliana} stipulates~\cite[Sec.~4.1.7]{lurie} that any displacement field $\neu$ expressed in the form
\begin{equation}
  \neu = 2(1-\nu)\Delta\bfg - \nabla(\operatorname{div}\bfg),\label{u:Gal}
\end{equation}
satisfies the Navier elastostatic system~\eqref{PDE:elas} with given body force density $\mu\nef$ provided the \emph{Galerkin vector potential} $\bfg$ satisfies the inhomogeneous biharmonic equation
\begin{equation}
  \Delta^2\bfg = \frac{\nef}{2(1-\nu)}\quad\text{in}\quad\R^d.  \label{g:biharm}
\end{equation}

A particular polynomial solution $\bfg\in[\Pcal_{N+4}]^d$ of~\eqref{g:biharm} can then be obtained for any given $\nef\in[\Pcal_N]^d$ by solving the vectorial bilaplacian equation \eqref{g:biharm} component-wise applying either of the approaches presented above in Sec.~\ref{sec:bilaplacian}, whereupon plugging $\bol g$ into~\eqref{u:Gal} yields a particular elastostatic displacement $\bol u$ generated by the force density $\mu\nef$.

\subsection{Stokes flows} As it turns out, the Galerkin vector potential representation~\eqref{u:Gal} applies also to the stationary inhomogeneous Stokes equation~\cite{kanwal:71}, which is formally identical to that of incompressible isotropic elasticity in the case $\nu=1/2$~\cite[Sec.\ 2.2.4]{Hsiao}. Indeed, a steady velocity field $\neu\in[\Pcal_{N+2}]^d$ and pressure field $p\in\Pcal_{N+1}$ solving
\begin{equation}
  \mu\Delta\neu - \nabla p = \nef , \qquad \operatorname{div}\neu=0\quad\text{in}\quad \R^d, \label{stokes}
\end{equation}
(where $\mu$ is here the dynamic viscosity of the fluid material) for a given body force density $\nef\in[\Pcal_N]^d$, can be expressed as
\begin{equation}
  \neu = \Delta\bfg - \nabla(\operatorname{div}\bfg), \qquad p = -\mu\Delta(\operatorname{div}\bfg), \label{u:Gal:stokes}
\end{equation}
where the {Galerkin vector} $\bfg\in[\Pcal_{N+4}]^d$ solves
\begin{equation}
  \Delta^2\bfg = \frac{\nef}{\mu}\quad\text{in}\quad \R^d. \label{g:biharm:stokes}
\end{equation}
In particular, the representation~\eqref{u:Gal:stokes} automatically satisfies the incompressibility constraint. Then, particular (polynomial) solutions can again be obtained for $(\neu,p)$ by solving (e.g. componentwise) the bilaplacian equation~\eqref{g:biharm:stokes} for $\bfg$ and using that solution in~\eqref{u:Gal:stokes}.

\subsection{Linearized Navier-Stokes equations}\label{sec:NStokes}
A common implicit-explicit time discretization of the incompressible Navier-Stokes equations in which the non-linear term is treated explicitly results in the problem: for a given $\nef\in[\Pcal_N]^d$, find a velocity field $\neu\in[\Pcal_{N}]^d$ and pressure field $p\in\Pcal_{N+1}$ solving
\begin{equation}
  (\Delta - \alpha^2)\neu - \mathrm{Re} \nabla p = \nef,\qquad \operatorname{div}\neu=0\quad\text{in}\quad \R^d,\label{linearized_navier_stokes}
\end{equation}
with $\alpha^2 \in \mathbb{C}$ a nonzero constant and where $\mathrm{Re}$ is the Reynolds number. The equation is alternately known~\cite{brinkman1949calculation,af2020fast} as the Brinkman equation, the linearized Navier-Stokes equation, or the modified Stokes equation.

We seek a solution $(\neu, p)$ in the form of the ansatz
\begin{equation}
    \neu = (\Delta + \alpha^2)(\Delta\bfg - \nabla (\operatorname{div}\bfg)),\qquad p = -\frac{1}{\mathrm{Re}}(\Delta^2 - \alpha^4) (\operatorname{div}\bfg). \label{Nstokes:pot}
\end{equation}
where $\bfg$ is a vector polynomial potential to be determined. We observe that any $\neu$ so defined satisfies the incompressibility condition $\operatorname{div} \neu = 0$, while a computation for the momentum equation reveals
\begin{equation}
\begin{split}
(\Delta - \alpha^2) \neu - \mathrm{Re}\nabla p &= (\Delta^2 - \alpha^4)(\Delta\bfg - \nabla (\operatorname{div} \bfg))  + (\Delta^2 - \alpha^4)\nabla (\operatorname{div}\bfg)\\
&= (\Delta^3 - \alpha^4 \Delta)\bfg,
\end{split}
\end{equation}
so seeking a vector potential $\bfg\in[\Pcal_{N+2}]^d$ satisfying (component-wise) the PDE
\begin{equation}\label{eq:brinkman_potential_pde}
    (\Delta^3 - \alpha^4 \Delta)\bfg = \nef
\end{equation}
will yield the desired solution pair $(\neu, p)\in [\Pcal_{N}]^d\times \Pcal_{N+1}$. We are not aware of the potential representation~\eqref{Nstokes:pot} previously appearing.

Equation~\eqref{eq:brinkman_potential_pde} can, in turn, be solved using Theorem~\ref{thm:soln_expr}
with $\Bcal_2 = -\alpha^4 \Delta$ and $\Rcal_2 = \Delta^3$, for which we have $p = r = 6,\, q = 2$. This requires solving a vector Poisson equation for each iteration in~\eqref{eq:neumannseries_generic_repr} using (for instance) the method of Sec.~\ref{sec:Laplace}.

\subsection{Anisotropic Laplacian, scalar wave equation}\label{laplace:aniso}

We now focus on the operator $\Delta_A$ defined as $\Delta_A u=\operatorname{div}(\bol{A}\nabla u)=A_{ij}\partial_{ij}u$, where $\bol{A}$ is a symmetric invertible $d\times d$ matrix; in particular we have $\Delta_I=\Delta$. To look for polynomial solutions $u\in\Pcal_{N+2}$ of the anisotropic PDE

\begin{equation}
\Delta_A u = f\quad\text{in}\quad\R^d,\label{eq:AnLap}
\end{equation}
 with $f\in\Pcal_{N}$, we define the anisotropic squared length $r^2_A$ of $\ner\in\R^d$ by
\begin{equation}
  r_A^2 := \ner^T\bol{A}^{-1}\ner. \label{rA:def}
\end{equation}
Generalizing  formulae~\eqref{Delta:ids}, it can be shown that
\begin{equation}
  \nabla (r_A^{2\ell}) = 2\ell r_A^{2\ell-2}\bol{A}^{-1}\ner,\quad\mbox{and}\quad
  \Delta_A (r_A^{2\ell}) = 2\ell(2\ell+d-2) r_A^{2\ell-2}.
\end{equation}
 Therefore, letting $v=r_A^{2\ell+2}h_n$ with $h_n\in\Hcal_n$ for some $\ell\in\mathbb N_0$ and $p\in\mathbb N_0$, we have
\begin{equation}
  \Delta_A v = \gamma_{\ell}^{n}r_A^{2\ell}h_n\quad\text{in}\quad\R^d, \label{DeltaAP}
\end{equation}
with the coefficients $\gamma_\ell^n$ again defined in~\eqref{DeltaP}. The proof of~\eqref{DeltaAP} is essentially identical to that of~\eqref{DeltaP}, so it is omitted for conciseness.

The aforementioned properties then yield the solution
\begin{equation}
  u = \sum_{\ell=0}^m c_\ell r_A^{2\ell+2} \Delta_A^\ell f, \label{sol:AnLap}
\end{equation}
of~\eqref{eq:AnLap} for $f\in\Hcal_{n}$ with the coefficients $\{c_\ell\}_{\ell=0}^m$ again recursively defined by~\eqref{Delta:ck}. As in the case of the isotropic Poisson equation~\eqref{poisson:pb}, the general solution of~\eqref{eq:AnLap} for an arbitrary $f\in\Pcal_{N}$ can be obtained by expanding $f$ in a basis of homogeneous polynomials.

If the matrix $\bol{A}$ is symmetric positive definite, $\Delta_A$ is the anisotropic Laplacian describing anisotropic conductivity properties and~\eqref{eq:AnLap} the anisotropic Poisson equation. Alternatively, if $\bol{A}$ is sign-indefinite with $d-1$ eigenvalues of the same sign and the remaining eigenvalue of the opposite sign, the PDE format~\eqref{eq:AnLap} and associated polynomial solution~\eqref{sol:AnLap} pertain to a scalar wave equation (where the last $d$-th coordinate in the principal axes of $\bol{A}$ can be treated as the time, the remaining $d-1$ ones as space coordinates, and constitutive properties in space are allowed to be anisotropic). In this case, $r^2_A$ given by~\eqref{rA:def} evaluates (a possibly-anisotropic version of) the Minkowski space-time squared length.

\subsection{Anisotropic elastostatics}

Finally, we consider a general anisotropic elastic material, whose constitutive behavior is described by the 4th-order elasticity tensor $\bol{\Ccal}$. The Cartesian components $\Ccal_{ijk\ell}$ of $\bol{\Ccal}$ satisfy the usual major and minor symmetries $\Ccal_{ijk\ell}=\Ccal_{k\ell ij}=\Ccal_{jik\ell}$ ($1\leq i,j,k,\ell\leq d$, with $d=2$ or $3$). Any elastostatic displacement field $\neu$ generated in such medium by a given body force density $\nef$ satisfies
\begin{equation}
  -\Ccal_{ijk\ell}\partial_{j\ell}u_k = f_i\quad\text{in}\quad \R^d, \qquad  1\leq i\leq d, \label{elas:aniso}
\end{equation}
where Einstein's implicit summation convention on repeated indices is used. Following classical solution methods for anisotropic elasticity (see e.g. \cite[Chap.~1]{mura}), let the \emph{Christoffel matrix} $K_{ik}=K_{ik}(\bol{\xi})$ be defined for any $\bol{\xi}\in\R^d$ by $K_{ik}(\bol{\xi})=\Ccal_{ijk\ell}\xi_j\xi_{\ell}$, so that $K_{ik}(-\mathrm{i}\partial)$ is the differential operator on the left-hand side of~\eqref{elas:aniso}. 

The matrix $K_{ik}(\bol{\xi})$ is positive definite and, together with its determinant $E(\bol{\xi}):=\text{det}(N_{ik}(\bol{\xi}))$ and its adjugate matrix $N_{ik}(\bol{\xi})$, satisfies the general matrix identity~\cite[Chap.~4, Sec.~4]{strang:linear}
\begin{equation}
  K_{ik}(\bol{\xi})N_{kj}(\bol{\xi})=\delta_{ij}E(\bol{\xi}). \label{christoffel:inverse}
\end{equation}
Since $\bol{\xi}\mapsto K_{ik}(\bol{\xi})$ is here homogeneous with degree 2, the scalar function $\bol{\xi}\mapsto E(\bol{\xi})$ and the matrix-valued function $\bol{\xi}\mapsto N_{ik}(\bol{\xi})$ are moreover homogeneous with respective degrees $2d$ and $2d-2$.

Then, let a displacement field $\neu$ be sought in the form
\begin{equation}
  u_i = N_{ij}(-\mathrm{i}\partial)g_j \label{galerkin:aniso}
\end{equation}
in terms of a vector potential $\bfg$. By using the above ansatz in~\eqref{elas:aniso} written in the Fourier domain, we find that~\eqref{galerkin:aniso}) defines a solution of~\eqref{elas:aniso} provided  $\bfg$ satisfies the following componentwise-scalar differential equation of order $2d$:
\begin{equation}
  E(-\mathrm{i}\partial)\bfg = \nef \label{2dorder}
\end{equation}

In analogy with the Galerkin representation of isotropic elastostatics, particular solutions $\neu\in[\Pcal_{N+2}]^d$ of~\eqref{elas:aniso} for given (polynomial) right-hand sides $\nef\in[\Pcal_{N}]^d$ are therefore obtained by solving equation~\eqref{2dorder} for $\bfg$ and substituting $\bfg$ into~\eqref{galerkin:aniso}. For general anisotropic materials, this task amounts to solving $d$ inhomogeneous scalar PDEs of $2d$-th order instead of the inhomogeneous second-order $d$-dimensional vector PDE~\eqref{elas:aniso}. However, there are classes of anisotropic elastic materials for which an explicit factorization of $E(\bol{\xi})$ is known, in which case solving~\eqref{2dorder} reduces to sequentially solving lower-order scalar PDEs with polynomial right-hand sides. In particular, for three-dimensional transversely isotropic materials (characterized by five independent elastic constants), we have~\cite{wil}
\begin{equation}
  E(\bol{\xi}) = \prod_{i=1}^3 \big(\, A_i(\xi_1^2+\xi_2^2)+\xi_3^2 \,\big),
\end{equation}
where the (real, positive) constants $A_1,A_2,A_3$ are known in terms of the material elastic constants. The above factorization thus implies that equation~\eqref{2dorder} leads to sequentially solving three inhomogeneous anisotropic-Laplace equations (with three different tensors $\bol{A}$), a task to which the method of Sec.~\ref{laplace:aniso} applies.  Such factorizations are also available for other cases of elastic anisotropy, e.g. hexagonal crystals~\cite{kroner:53} and two-dimensional cubic crystals.

For the special case of isotropic elasticity, $N_{ij}(\bol{\xi})$ and $E(\bol{\xi})$ are found to be homogeneous with respective degrees 2 and 4 (for both $d=2$ and $d=3$), and their known closed-form expressions~\cite[Chap.~1]{mura} show that the foregoing method becomes identical to the Galerkin representation method used in Sec.~\ref{elas:iso}.

\section{Implementation details}\label{sec:numerical_aspects}

The methods proposed in the previous sections for constructing polynomial solutions have been implemented in the Julia~\cite{bezanson2017julia} library \href{https://github.com/IntegralEquations/ElementaryPDESolutions.jl}{\texttt{ElementaryPDESolutions.jl}}, made available under an MIT license on GitHub. The library is self-contained (i.e., it has no dependencies other than the Julia language itself), and consists of a few hundred lines of code; in particular, no symbolic computations are performed so that both the computation and the evaluation of the polynomial solutions are fully numerical procedures (i.e. the polynomial coefficients are regular numeric types such as single or double precision floating point numbers, long integers, etc).

In our implementation, polynomials are represented as a dictionary mapping a $d$-tuple $\alpha$ of exponents (the keys) to the corresponding coefficient $c_{\alpha}$ (the values). To allow for flexibility on both the spatial dimension and on the type of the coefficients, a generic \texttt{Polynomial\{N,T\}} type is defined, templated on both the ambient dimension \texttt{N} and on the numerical type of coefficient \texttt{T}; in Julia parlance, the \texttt{Polynomial} structure is said to be of parametric type. A bi-variate polynomial with double precision coefficients corresponds for instance to a \texttt{Polynomial\{2,Float64\}} object.

Given a polynomial \texttt{Q}, a polynomial solution \texttt{P} is obtained by means of the appropriate invocation \texttt{solve\_pde(Q,parameters)}, where \texttt{pde} corresponds to one of the supported partial differential equations (currently available choices for \texttt{pde} are \texttt{helmholtz}, \texttt{elastodynamics}, \texttt{maxwell}, \texttt{laplace}, \texttt{anisotropic\_laplace}, \texttt{anisotropic\_advect}, \texttt{anisotropic\_advect\_diffuse}, \texttt{bilaplace}, \texttt{elastostatics}, \texttt{stokes}, \texttt{brinkman}), and \texttt{parameters} are the numerical values of the physical parameters. For example, solving $\Delta P + 4 P = x^2y^3z$ is accomplished through the code in~\Cref{listing:helmholtz}.
\begin{listing}[t]
\begin{minted}[
    framesep=2mm,
    baselinestretch=1.2,
    frame=lines
]{julia-repl}
julia> f = Polynomial((2,3,1)=>1)
x²y³z
julia> u = solve_helmholtz(f,k=2)
0.375yz - 0.125y³z - 0.375x²yz + 0.25x²y³z
\end{minted}
\caption{Helmholtz with floating point coefficients}
\label{listing:helmholtz}
\end{listing}
For this example, the construction of \texttt{P} takes around a microsecond, and the returned polynomial can be evaluated at three-dimensional points in a few nanoseconds (on a 2022 MacBook Pro with a 2.3 GHz 8-core Intel Core i9 processor), making the library sufficiently fast for our main application of interest.

Because we avoid symbolic algebra, the coefficients of the polynomial solutions are subject to truncation errors if the intermediate stages of the computation cannot be represented exactly. While this is not necessarily a problem in itself (we did not observe catastrophic accumulation of truncation errors in the tests we performed using double precision floating point numbers), it may be convenient to either avoid truncation error altogether, or to obtain rigorous interval bounds on the coefficients which are computed. Both means of providing guarantees on the correctness of the computed coefficients are naturally supported by the \texttt{ElementaryPDESolutions.jl} library.

To illustrate how rational numbers can be used instead of floating point numbers, we consider again~\Cref{listing:helmholtz}, but modify $k$ to be of \texttt{Rational} type (Julia provides native support for rational numbers). The code snippet shown in~\Cref{listing:helmholtz-rational} illustrates how this is accomplished.
\begin{listing}[b]
\begin{minted}[
    framesep=2mm,
    baselinestretch=1.2,
    frame=lines
]{julia-repl}
julia> f = Polynomial((2,3,1)=>1)
x²y³z
julia> u = solve_helmholtz(f,k=Rational(2))
3//8yz - 1//8y³z - 3//8x²yz + 1//4x²y³z
\end{minted}
\caption{Helmholtz with rational coefficients. The double-slash notation \texttt{a//b} denotes the fraction $\frac{\texttt{a}}{\texttt{b}}$.}
\label{listing:helmholtz-rational}
\end{listing}
Note that because~\eqref{P:sol:helm} involves only iterated Laplacians of \texttt{Q} and divisions by powers of $k^2$, the intermediate stages of the computation involve only rational coefficients provided both $k$ and the coefficients of $Q$ are rational\footnote{For large polynomial orders, the $64$-bit integer types used by default in Julia may overflow. Since Julia provides support to multiple precision arithmetic by wrapping the GNU MP library~\cite{granlund2010gnu}, a simple fix is to use the \texttt{BigInt} type if needed (e.g. use \texttt{k=Rational\{BigInt\}(2)} instead of \texttt{k=Rational(2)}).}. The same is true for the other Helmholtz-like problems presented in this paper.
\begin{listing}[b]
\begin{minted}[
    framesep=2mm,
    baselinestretch=1.2,
    frame=lines
]{julia-repl}
julia> using IntervalArithmetic
julia> f = Polynomial((2,3,1)=>1)
x²y³z
julia> u = solve_helmholtz(f,k=Interval(pi))
[0.0249638, 0.0249639]yz - [0.0205319, 0.020532]y³z -
[0.0615958, 0.0615959]x²yz + [0.101321, 0.101322]x²y³z
\end{minted}
\caption{Helmholtz with interval coefficients. The square brackets in the coefficients of \texttt{P} provide a lower and upper bound on the value.}
\label{listing:helmholtz-interval}
\end{listing}

When the problem parameters are not rational numbers, we may either approximate them by rationals to the desired precision and proceed as shown before, or use interval arithmetic~\cite{moore2009introduction} to propagate error bounds on floating point operations. Using the \texttt{IntervalArithmetic.jl} package, we can easily construct an \texttt{Interval} representation for the problem parameters (e.g. \texttt{k = Interval(pi)}), and then pass it to our solver as shown in~\Cref{listing:helmholtz-interval}. Interestingly, due to the generic nature of our \texttt{Polynomial} type, the code works as is even when the coefficients of the polynomials are \texttt{Interval} objects. The computed solution has coefficients which are intervals instead of numbers, thus providing error bounds on the coefficient values due to truncation.\enlargethispage*{3ex}

Although we have used the Helmholtz equation to illustrate some of the functionality and possible pitfalls, the library supports many other PDEs under a very similar API. An example illustrating how one can solve the (vectorial) Stokes system, where the solution is much less trivial to obtain manually, is shown in~\Cref{listing:stokes-rational}.
\begin{listing}[t]
\begin{minted}[
    escapeinside=||,
    framesep=2mm,
    baselinestretch=1.2,
    frame=lines
]{julia-repl}
julia> f = (Polynomial((1,1)=>1),Polynomial((1,0)=>1))
(xy, x)
julia> u,p = solve_stokes(f,μ=Rational(2));
julia> u
 -1//96y³ - 1//32x²y + 1//64x³y + 5//192xy³
 1//32xy² + 5//96x³ - 5//768x⁴ - 3//128x²y² - 5//768y⁴
julia> p
-1//4xy - 1//12y³ - 1//4x²y
\end{minted}
\caption{Stokes with rational coefficients}
\label{listing:stokes-rational}
\end{listing}

Examples of user-specified anisotropy in Laplace, advection, and advection-diffusion equations are given in~\Cref{listing:anisotropy}.
\begin{listing}[t]
\begin{minted}[
    escapeinside=||,
    framesep=2mm,
    baselinestretch=1.2,
    frame=lines
]{julia-repl}
julia> using StaticArrays
julia> A = SMatrix{2,2,Rational{Int64}}(2//1, 1//1, 1//1, -3//1);
julia> β = SVector{2,Rational{Int64}}(2//1, 1//1);
julia> g = Polynomial([(1, 1) => 2//1])
2//1xy
julia> v = solve_anisotropic_laplace(A, g)
-3//784x⁴ + 13//196x³y - 13//294xy³ + 5//98x²y² - 1//588y⁴
julia> v = solve_anisotropic_advect(β, g)
13//125xy² - 26//375y³ + 56//125x²y - 28//375x³
julia> v = solve_anisotropic_advect_diffuse(A, β, g)
-1424//3125y - 2848//3125x - 698//625xy - 382//625y² + 132//625x² + 13//125xy²
- 26//375y³ + 56//125x²y - 28//375x³
\end{minted}
\caption{Anisotropic models. The \texttt{StaticArrays} package includes routines for sufficiently small dimension-$d$ which explicitly invert the (note, additionally, non-positive) rational matrix \texttt{A} as needed, in a manner compatible with \texttt{Rational} types. (\texttt{StaticArrays} is \emph{not} a dependency of the package but \texttt{using StaticArrays} will lead to faster construction and enable \texttt{Rational} support for these problems.)}
\label{listing:anisotropy}
\end{listing}

Finally, in~\Cref{listing:maxwell} we show an example for the (three-dimensional) Maxwell system, where the coefficients of the solution are complex numbers.
\begin{listing}[t]
\begin{minted}[
    escapeinside=||,
    framesep=2mm,
    baselinestretch=1.2,
    frame=lines
]{julia-repl}
julia> J = (Polynomial((2,1,0)=>1),Polynomial((1,0,0)=>1),Polynomial((0,0,0)=>1))
(x²y, x, 1)
julia> E,H = solve_maxwell(J,μ=2);
julia> E
((-0.0 - 1.0im)x²y, (-0.0 - 2.0im)x, (-0.0 - 1.0im))
julia> H
(0, 0, -1.0 + 0.5x²)
\end{minted}
\caption{Maxwell system. Additional keyword arguments for $\epsilon$ and $\omega$ can be passed to \texttt{solve\_maxwell}; by default their value is one.}
\label{listing:maxwell}
\end{listing}

\section{Conclusions}

We presented a general methodology for finding polynomial solutions to various
linear, constant coefficient PDEs, in the presence of a polynomial source term.
The method, based on a formal Neumann series of the differential operator, significantly generalizes related past techniques used
in specific cases like the Helmholtz-type problems~\cite{Dangal2017} and allows for treatment of quite general PDOs, requiring an inverse operator only for the lowest-order part of the PDO. Furthermore, to the best of our
knowledge, we propose the first general method for obtaining polynomial
particular solutions of equations incorporating a
divergence constraint, such as the Stokes, linearized Navier-Stokes and and Maxwell systems. We expect that
the presented methods and the accompanying Julia library will be applicable to
other PDE models not treated here, such as the poroelasticity system, and prove
useful to others for developing PDE solution methods that require particular
polynomial solutions.

\section{Acknowledgements}
The authors are grateful for the detailed and insightful comments by one reviewer, which much helped them in expanding, clarifying and overall significantly improving this article.

\bibliographystyle{plainurl}
\bibliography{References}
\end{document}